\documentclass[12pt,a4paper]{amsart}
\usepackage{amssymb,amsmath}
\usepackage{graphicx} 

 \usepackage{colortbl}
  \usepackage{cancel}
\numberwithin{equation}{section}
  
\definecolor{red}{rgb}{1,0,0}

\definecolor{blue}{rgb}{0,0,1}

\def\rife#1{(\ref{#1})}
\newtheorem{example}{Example}[section]

\newtheorem{lemma}[example]{\sc Lemma}

\newtheorem{prop}[example]{\sc Proposition}
\newtheorem{theo}[example]{\sc Theorem}

\def\proof{\noindent{\sc  Proof.}\quad}
\def\sqr#1#2{{\vcenter{\vbox{\hrule height.#2pt 
													\hbox{\vrule width.#2pt height#1pt \kern#1pt
													\vrule width.#2pt}
 \hrule height.#2pt}}}}
\def\finedim{{\unskip\nobreak\hfil\penalty50
\hskip2em\hbox{}\nobreak\hfil\hbox{$ \sqr44$ \qquad}
\parfillskip=0pt \finalhyphendemerits=0\par\medskip}}
\def\disp{\displaystyle}
\def\R{I \!\!R}
\def\N{I \!\! N}
\def\elle#1{L^{#1}(\Omega)}

\def\w{W^{1,2}_0(\Omega)}
\def\ww{W^{1,1}_0(\Omega)}

\def\io{\int_{\Omega}}
\def\ik{\int_{\{k\leq|\un|\}}}
    
\def\ie{\int_{E}}
\def\norma#1#2{\|#1\|_{\lower 4pt \hbox{$ \scriptstyle #2$ }}}
\def\div{{\rm div}}

\def\un{u_{n}}
\def\zn{z_{n}}

\def\fn{f_{n}}
\def\Fn{F_{n}}

\def\deb{\rightharpoonup}

\begin{document}
\title[\sc An elliptic system with degenerate coercivity]
{\sc  An elliptic system with degenerate coercivity}
\author[\sc L. Boccardo, G. Croce, C. Tanteri]
{\sc  Lucio Boccardo, Gisella Croce, Chiara Tanteri}
\address{La sapienza Universit\`a di Roma.}
\email{boccardo@mat.uniroma1.it}
\address{Universit\'e du Havre}
\email{gisella.croce@univ-lehavre.fr}
\address{\'Ecole polytechnique f\'ed\'erale de Lausanne}
\email{chiara.tanteri@epfl.ch}
\maketitle
 
 \rightline{\it a  Bernard, nostro maestro}
 \rightline{\it 
 \footnote{ (see \cite{d1}, \cite{d2}, \cite{bd84}, \cite{bd89},\cite{cd},  \cite{dt1}, \cite{dt2},) }
 }
 
  

\section{Introduction}
\subsection{\sc Setting}
In this paper we study the existence of solutions of the  {degererate} elliptic system
\begin{equation}
\label{sist}
\left\{
\begin{array}{c}
 -\div\bigg(\dfrac{\,a(x) \nabla u }{(b(x)+|z|)^2}\bigg)+u=f(x),
 \\  
 \\
 -\div\bigg(\dfrac{\, A(x) \nabla z  }{(B(x)+|u|)^2}\bigg)+z=F(x),
\end{array}
\right.
\end{equation}
where
 $\Omega$ is a bounded, open subset of $\R^{N}$, with $N > 2$,
   $\,a(x) $  and $ \,A(x) $ are  measurable matrices such that, for 
   $\alpha,\,\beta\in\R^+$,
   \begin{equation}
\label{a}
 \alpha|\xi|^2\leq \,a(x)\xi\xi,\;
 \alpha|\xi|^2\leq  \,A(x)\xi\xi;\quad 
 |\, a(x)|\leq\beta ,\;
 |\, A(x)|\leq\beta .
\end{equation}
Moreover we assume
\begin{equation}
\label{b}
0<\lambda\leq b(x),\;B(x)\leq\gamma,
\end{equation} 
for some 
   $\lambda,\,\gamma\in\R^+$ and
\begin{equation}
\label{f}
f(x),\,F(x) \in\elle 2 .
\end{equation}
\begin{theo}\label{theo}
Under the assumptions \rife{a}, \rife{b}, \rife{f}, 
there exist  $u\in\ww$ and $z\in\ww$, distributional solutions of the system \rife{sist}.
\end{theo}
\subsection{\sc Comments}
First of all, we note that existence of solutions belonging to the nonreflexive space $\ww$ is not so usual in the study of elliptic problems.
Recently the existence of solutions in $\ww$ was proved in
\cite{bco3}, \cite{bco1}, \cite{bco2}, for elliptic scalar problems with degenerate coercivity (so that this paper is an extension to the systems of some of those results) and in some borderline cases of the Calderon-Zygmund theory of nonlinear Dirichlet problems in \cite{bg2012}.

The main difficulty of the problem is that
even if the differential operator is well defined
between $\w$ and its dual, it is not coercive on $\w$: degenerate coercivity means that  when $|v|$ is ``large'', 
  $\frac{1}{(b(x)+|v|)^2}$ goes to zero: for  an   explicit example see \cite{Po}.
  
 The study of problems involving degenerate equations
begins with the paper \cite{bdo} 
 and it is  developed in \cite{bb}, \cite{croce1}, \cite{croce2}, \cite{croce3}, \cite{bco3}, \cite{bco1}, \cite{bco2}
 (see also  \cite{bc-libro}) 
  
 \section{Existence}
\subsection{\sc A priori estimates}  
The first existence result is concerned with the case of a bounded data.
  
We recall the following definitions.
$$
T_k(s)= 
\left\{
\begin{array}{rcl}
s,
& \mbox{if $|s|\leq k$;} \\
k\frac{s}{|s|},
& \mbox{if $|s|> k$;} \\
\end{array}
\right. 
\qquad
G_k(s)=s-T_k(s).
$$
\begin{prop}\label{bar}\rm  
Let $\rho>0$, $\sigma>0$ and $g,\,G\in\elle\infty$.
Then there exist    weak solutions $w$, $W$ belonging to $\w$ of the system
$$
 \left\{
\begin{array}{c}
w\in\w\cap\elle\infty:\;
  -\div\bigg(\dfrac{\, a(x) \nabla w}{(b(x)+|T_\rho(W)|)^2}\bigg)+w=g(x),
 \\  
 \\
W\in\w\cap\elle\infty:\;-\div\bigg(\dfrac{\, A(x) \nabla W}{(B(x)+|T_\sigma(w)|)^2}\bigg)+W=G(x). 
\end{array}
\right.
$$
\end{prop}
\proof
The  existence   
is a consequence of the Leray-Lions theorem (or Schauder theorem) since the principal part is not degenerate, thanks to the presence of $T_\rho$ and $T_\sigma$.
Moreover, if we take 
$G_h(w)$ as test function in the first equation and 
$G_k(W)$ as test function in the second equation, we have, dropping two positive terms,
$$
 \left\{
\begin{array}{c}
\disp
\io[|w|-|g(x)|]|G_h(w)|\leq0,
\\
\\
\disp
\io[|W|-|G(x)|]|G_k(w)|\leq0.
\end{array}
\right.
$$
Then the choice $h=\norma{g}{\elle\infty}$, $k=\norma{G}{\elle\infty}$ implies
$$
 \left\{
\begin{array}{c}
\disp
|{w}|\leq \norma{g}{\elle\infty},
\\
\disp
|{W}|\leq \norma{G}{\elle\infty}.
\end{array}
\right.
$$
Thus, if we set $\rho=\norma{g}{\elle\infty}$ and
$\sigma=\norma{G}{\elle\infty}$, we can say that $w$ and $W$ are   bounded weak solutions of the system
$$
 \left\{
\begin{array}{c}
w\in\w\cap\elle\infty:\;
  -\div\bigg(\dfrac{\, a(x) \nabla w}{(b(x)+|W|)^2}\bigg)+w=g(x),
 \\  
 \\
W\in\w\cap\elle\infty:\;-\div\bigg(\dfrac{\, A(x) \nabla W}{(B(x)+|w|)^2}\bigg)+W=G(x). 
\end{array}
\right.
$$
\finedim
Now we define
$$
\fn=\frac{f}{1+\frac1n|f|},\qquad
\Fn=\frac{F}{1+\frac1n|F|},
$$
so that  
 \begin{equation}
\label{banga}
\norma{\fn-f}{\elle2}\to0,\qquad
\norma{\Fn-F}{\elle2}\to0.
\end{equation}
Thanks to the   Proposition \ref{bar}, there 
  exists a solution $(\un,\zn)$ of the system
\begin{equation}
\label{sistn}
\left\{
\begin{array}{c}
\un\in\w:\;
  -\div\bigg(\dfrac{\, a(x) \nabla\un}{(b(x)+|\zn|)^2\,}\bigg)+\un=\fn(x),
 \\  
 \\
\zn\in\w:\;-\div\bigg(\dfrac{\, A(x) \nabla\zn}{(B(x)+|\un|)^2}\bigg)+\zn=\Fn(x), 
\end{array}
\right.
\end{equation}
Now we prove our first estimates.
\begin{lemma}
The sequences $\{\un\}$ and $\{\zn\}$ are bounded in $\elle2$.
\end{lemma}
\proof
We take   $G_k(\un)$ as a test function in the first equation and we have
\begin{equation}
\label{visto}
\alpha\disp
\io\dfrac{|\nabla G_k(\un)|^{2}}{( b(x)+|\zn|)^{2}}
+
\io |G_k(\un)|^{2}
\leq
\io |f|\,|G_k(\un)|
\end{equation}
If we drop the first positive term and we use the H\"older inequality, then we have
\begin{equation}
\label{india}
\bigg[\io |G_k(\un)|^{2}\bigg]^\frac12
\leq
\bigg[ \ik|f|^2\bigg]^\frac12.
\end{equation}
Similar estimates hold true for  $\zn$. 
In particular, taking $k=0$, we have the boundedness of the sequences $\{\un\}$ and $\{\zn\}$ in $\elle2$. So we have
that there exist $u$, $z$ such that,
up to subsequences, 
 \begin{equation}
\label{aula}
\un\deb\,u,\quad
\zn\deb\,z \qquad
\hbox{weakly in $\elle2$}.
\end{equation}
Then if we drop the   second term in \rife{visto}, we have
\begin{equation}
\label{bangalore}
\alpha\disp
\io\dfrac{|\nabla G_k(\un)|^{2}}{( b(x)+|\zn|)^{2}}
\leq
\ik|f|^2.
\end{equation}
A similar estimate  for  $\zn$ comes from the second equation. 
\finedim
\begin{lemma}\label{andrea}
The sequences $\{\un\}$ and $\{\zn\}$ are bounded in $\ww$.
\end{lemma}
\proof
A consequence of \rife{bangalore} and of   the H\"older inequality is
$$
\begin{array}{c}
\disp
\io|\nabla G_k(\un)|
=\io\dfrac{|  \nabla G_k(\un)| }{( b(x)+|\zn|) }\,( b(x)+|\zn|)
\\
\disp
\leq
\bigg[\ik\dfrac{|f|^2}\alpha\bigg]^\frac12
\big(\norma{b}{\elle2}+\norma{f}{\elle2} \big).
\end{array}
$$
Similar estimates hold true for  $\zn$.  
In particular, with $k=0$, we have
\begin{equation}
\label{bolt}
\begin{array}{c}
\disp
\io|\nabla\un|\leq
 \dfrac{\norma{f}{\elle2}\,\big(\norma{b}{\elle2}+\norma{f}{\elle2} \big)}{\alpha^\frac12} ,         
         \\
   \\      
   \disp     
\io|\nabla\zn|\leq
 \dfrac{\norma{F}{\elle2}\,\big(\norma{b}{\elle2}+\norma{f}{\elle2} \big)}{\alpha^\frac12}  .
\end{array}
\end{equation}
\finedim
Now we improve the convergence \rife{aula}.
\begin{lemma}\label{mennea}
The sequences $\{\un\}$ and $\{\zn\}$ are compact in $\elle2$.
\end{lemma}
\proof
The estimates \rife{bolt} imply, thanks to the Rellich embedding Theorem, the $L^1$ compactenss and then the a.e. convergences
\begin{equation}
\label{lewis}
\un(x)\to\,u(x),\qquad
\zn(x)\to\,z(x).
\end{equation}
Now we use the Vitali Theorem: since we have the a.e. convergences \rife{lewis}, the compactness is achieved if we prove the equiintegrability.

Let    $E$ be a measurable subset of $\Omega$. 
Since $\un=T_k(\un)+G_k(\un)$, we have (we use \rife{india})
$$
\begin{array}{c}
 \disp
\ie |\un|^2
\leq 2\ie|T_k(\un)|^2+2\ie|G_k(\un)|^2
\\
\disp
\leq 2\,k^2\,|E| +2\io|G_k(\un)|^2
\leq
\\
\disp
2\,k^2\,|E| +2\ik|f|^2,
\end{array}
$$ 
where $|E|$ denotes the measure of $E$.
Now we recall that a consequence of Lemma \ref{andrea} is that the sequence $\{\un\}$ is bounded in $\elle1$, so that if we fix $\epsilon>0$, there exists $k_\epsilon$ such that
(uniformly with respect to $n$)
$$
\ik|f|^2\leq\epsilon, \quad
k\geq\,k_\epsilon.
$$
Then
$$
\ie |\un|^2
\leq
2\,k^2\,|E| +2\epsilon
$$
implies
$$
\lim_{|E|\to0}
\ie |\un|^2
\leq
2\epsilon, \;\hbox{uniformly with respect to $n$}.
$$
Similar inequality holds true for  $\zn$. 
\finedim
\begin{lemma} \label{berruti}
The sequences $\{\un\}$ and $\{\zn\}$ are weakly compact in $\ww$.
\end{lemma}
\proof
Here we follow \cite{bco1}, \cite{bco2}.
Let  again $E$ be a measurable subset of $\Omega$, and let $i$ be in $\{1,\ \ldots,\ N\}$. Then
$$
\disp
\int_{E} |\partial_{i}\un|
 \leq 
\disp 
\int_{E} |\nabla \un|
=
\int_{E} \frac{|\nabla \un|}{b(x)+|\zn|}\,(b(x)+|\zn|)
$$
$$
\leq
\disp \bigg[ 
\io\frac{|\nabla \un|^2}{(b(x)+|\zn|)^{2}}
\bigg]^\frac 12
\bigg[\int_{E} (b(x)+|\zn|)^{2}\bigg]^\frac{1}{2}
$$
$$
\leq
\disp
\bigg[\frac{1}{\alpha}\io |f|^2\bigg]^{\frac12}
\bigg\{
\bigg[\int_{E} b(x)\bigg]^\frac{1}{2}
+ \bigg[\int_{E}  |\zn|^{2}\bigg]^\frac{1}{2} \bigg\}\,,
$$
where we have used the inequality \rife{bangalore} in the last passage. Since the sequence $\{\un\}$ is compact in $\elle{2}$, we have that the sequence $\{\partial_{i}\un\}$ is equiintegrable. Thus, by Dunford-Pettis theorem, and up to subsequences, there exists $Y_{i}$ in $\elle1$ such that 
$\partial_{i} \un$ weakly converges to $Y_{i}$ in $\elle1$. Since $\partial_{i}\un$ is the distributional derivative of $\un$, we have, for every $n$ in $\N$,
$$
\io \partial_{i} \un\,\phi = -\io \un\,\partial_{i} \phi\,,
\quad
\forall\; \phi \in C^{\infty}_{0}(\Omega)\,. 
$$
We now pass to the limit in the above identities, using that $\partial_{i}\un$ weakly converges to $Y_{i}$ in $\elle1$, and that $\un$ strongly converges to $u$ in $\elle2$; we obtain
$$
\io Y_{i}\,\phi = -\io u\,\partial_{i} \phi\,,
\quad
\forall \; \phi \in C^{\infty}_{0}(\Omega)\,, 
$$
which implies that $Y_{i} = \partial_{i} u$, and this result is true for every $i$. Since $Y_{i}$ belongs to $\elle1$ for every $i$, $u$ belongs to $W^{1,1}_{0}(\Omega)$.
A similar result holds true for  $\zn$. 
\finedim
Thus, thanks to Lemma \ref{mennea}  and Lemma \ref{berruti},
we can improve the convergence \rife{aula}:
\begin{equation}
\label{owens}
\left\{
\begin{array}{c}
\un \hbox{ converges weakly in } \ww 
\hbox{ and strongly in } \elle2 \hbox{ to } u,
      \\
     \zn \hbox{ converges weakly in } \ww 
\hbox{ and strongly in } \elle2 \hbox{ to } z.   
\end{array}
\right.
\end{equation}
\finedim
\subsection{\sc Proof of Theorem \ref{theo} - } 
First of all, we use the equiintegrability proved in Lemma \ref{berruti}: fix $\varepsilon>0$, there exists 
$\delta(\varepsilon)>0$ such that, for every measurable subset $E$ with $|E|\leq\delta(\varepsilon)$, we have
$$
\ie|\nabla\un|\leq\varepsilon.
$$
Taking into account the absolute continuty of the Lebesgue integral, we have, for some $\tilde\delta(\varepsilon)>0$,
$$
\ie|\nabla\un|\leq\varepsilon,\quad
\ie|\nabla u|\leq\varepsilon ,
$$
for every measurable subset $E$ with $|E|\leq\tilde\delta(\varepsilon)$.

On the other hand, since $|\Omega|$ is finite and the sequence  
$$
D_n=\dfrac{\, a(x) }{(b(x)+|\zn|)^2\,}
$$
converges almost everywhere (recall   \rife{owens}),
 the Egorov theorem says that 
for every $q> 0$, there exists a measurable subset $F$ of $\Omega$ such that $|F|< q$ , and $D_n$ converges to $D$ uniformly on $\Omega\setminus F$.
We choose $q=\tilde\delta$ so that we have, for every $\varphi\in{\rm Lip}(\Omega)$,
$$
\begin{array}{c}
     \disp
\bigg|\io[D_n \nabla\un\nabla\varphi
-D \nabla u\nabla\varphi]\bigg|
\\
\disp
\leq
\bigg|\int_{\Omega\setminus F}[D_n \nabla\un\nabla\varphi-D \nabla u\nabla\varphi]\bigg|
+\bigg|\int_{F}[D_n \nabla\un\nabla\varphi-D \nabla u\nabla\varphi]\bigg|
\\
\disp
\leq
\Big|\int_{\Omega\setminus F}[D_n \nabla\un\nabla\varphi-D \nabla u\nabla\varphi]\Big|
+ \frac{\beta}{\lambda^2}
\norma{\,|\nabla\varphi|\,}{\elle\infty}
\Big[  \int_{F}|\nabla\un| + \int_{F}| \nabla u ]\Big]
\\
\disp
\leq
\Big|\int_{\Omega\setminus F}[D_n \nabla\un\nabla\varphi-D \nabla u\nabla\varphi]\Big|
+ 2\varepsilon\frac{\beta}{\lambda^2}
\norma{\,|\nabla\varphi|\,}{\elle\infty},
\end{array}
$$
which proves that
 \begin{equation}
\label{lore}
\io\,\dfrac{\,a(x)\,\nabla\un\nabla\varphi}{(b(x)+|\zn|)^2\,}
\;\to\;\io\,
\dfrac{\,a(x)\,\nabla u\nabla\varphi}{(b(x)+|z|)^2} .
\end{equation}
Thus,
thanks to the above limit, \rife{banga} and Lemma \ref{mennea},
 it is possible to pass to the limit in the weak formulation of \rife{sistn}, for every $\varphi,\;\psi \in  {\rm Lip}(\Omega)$,
\begin{equation}
\label{chiarsella}
\left\{
\begin{array}{c}
\disp
\io\dfrac{\, a(x) \nabla\un \nabla \varphi}
{(b(x)+|\zn|)^2\,}
+\io\un\,\varphi=\io\fn(x)\,\varphi,
 \\  
 \disp
\io \dfrac{\, A(x) \nabla\zn  \nabla \psi}{(B(x)+|\un|)^2} 
+\io\zn\, \psi=\io\Fn(x) ;
\end{array}
\right.
\end{equation}
and we prove that $u$ and $z$ are solutions of our system, in the following distributional sense
 \begin{equation}
\label{gennaio}
\left\{
\begin{array}{c}
\disp
\io\dfrac{\, a(x) \nabla u \nabla \varphi}{(b(x)+|z|)^2}
+\io u\,\varphi=\io f(x)\,\varphi,
\quad\forall\;\varphi\in\, {\rm Lip}(\Omega);
 \\  
 \disp
\io \dfrac{\, A(x) \nabla z  \nabla\psi}{(B(x)+|u|)^2} 
+\io z\, \psi=\io F(x)\,\psi, 
\quad\forall\;\psi\in\, {\rm Lip}(\Omega).
\end{array}
\right.
\end{equation}
\finedim
Now we show that, in the above   definition of solution, it is possible to use less regular test functions: it possible to use functions only belonging to $\w$. 
\begin{prop}
The above functions $u$ and $z$ are solutions of our system, in the following sense
 \begin{equation}
\label{25}
\left\{
\begin{array}{c}
\disp
\io\dfrac{\, a(x) \nabla u \nabla v}{(b(x)+|z|)^2}
+\io u\,v=\io f(x)\,v,
\quad\forall\;v\in\w;
 \\  
 \disp
\io \dfrac{\, A(x) \nabla z  \nabla w}{(B(x)+|u|)^2} 
+\io z\, w=\io F(x)\,w, 
\quad\forall\;w\in\w.
\end{array}
\right.
\end{equation}
\end{prop}
\proof
In order to avoid technicalities, here we also assume that
\begin{equation}
\label{bernard}
\hbox{$a(x)$ and $A(x)$ are scalar functions.}
\end{equation}  
We start with the following inequalities (we use \rife{bangalore} with $k=0$)
$$
\io\bigg|\dfrac{\,a(x)\nabla\un }{(b(x)+|\zn|)^2\,}\bigg|^2
\leq
\frac{\alpha^2}{\lambda^2}
\io\dfrac{\,|\nabla\un|^2\,}{(b(x)+|\zn|)^2\,} 
\leq
\frac{\alpha^2}{\lambda^2}
\io|f|^2.
$$
Thus, up to subsequences, we can say that, for some $\Psi\in(\elle2)^N$,
 \begin{equation}
\label{dacorogna}
\io\dfrac{\,a(x)\nabla\un }{(b(x)+|\zn|)^2\,}\,\Phi
\;\to\;
\io\Psi\,\Phi,
\end{equation}
for every $\Phi\in(\elle2)^N$.
Now we compare the limit  \rife{lore} with the limit  \rife{dacorogna}, taking $\Phi=\nabla\varphi$, and we deduce that
$$
\io\bigg[\dfrac{\,a(x)\,\nabla u\, }{(b(x)+|z|)^2}-\Psi\bigg]\, \Phi=0.
$$
Thus we proved that
$$
\dfrac{\,a(x)\nabla\un }{(b(x)+|\zn|)^2\,}
\;\hbox{ weakly converges in $(\elle2)^N$ to }\;
\dfrac{\,a(x)\nabla u }{(b(x)+|z|)^2\,},
$$
which allows us to pass to the limit in \rife{chiarsella}
only assuming $\varphi,\;\psi \in\w$.
\section*{Acknoledgments}
This paper contains  the unpublished part  of the results presented  by the first author 
in a talk at the conference ``Calculus of Variations and Differential Equations - Conf\'erence en l'honneur du 60\`eme anniversaire de Bernard Dacorogna" (Lausanne, 10-12 juin 2013).

 \end{document}